\documentstyle[amssymb,11pt]{article}

\title{A saturated model of
an unsuperstable theory of cardinality greater than its theory has the small 
index property}

\author{Garvin
Melles\thanks{Would like to thank Ehud Hrushovski
for supporting him with funds from NSF Grant DMS 8959511} \\Hebrew
University of Jerusalem \and Saharon Shelah\thanks{partially supported
by the U.S.-Israel Binational 
Science Foundation. Publ. 450}\\Hebrew
University of Jerusalem\\Rutgers University}

\newtheorem{theorem}{Theorem}[section]
\newtheorem{defi}[theorem]{Definition}
\newtheorem{lemma}[theorem]{Lemma}
\newtheorem{coro}[theorem]{Corollary}
\newtheorem{notation}[theorem]{Notation}

\newcommand{\sub}{\subseteq}
\newcommand{\proof}{{\sc proof} \hspace{0.1in}}
\newcommand{\iopp}{{\mathop{\bigcup}}{\!\!\!\!{|}}}
\newcommand{\siopp}{{\mathop{\bigcup}}{\!\!\!{|}}}
\newcommand{\niopp}{\not\!\!\iopp}
\newcommand{\noindep}[1]{\mathop{\niopp}\limits_{\textstyle{#1}}\ }
\newcommand{\indep}[1]{\mathop{\iopp}\limits_{\textstyle{#1}}\ }
\newcommand{\sindep}[1]{\mathop{\siopp}\limits_{\textstyle{#1}}\ }

\begin{document}
\mathsurround=.1cm
\maketitle

\begin{abstract}
A model $M$ of cardinality $\lambda$ is said to have the small index
property if for every $G\sub Aut(M)$ such that $[Aut(M):G]\leq\lambda$
there is an $A\sub M$ with $\big|A\big|<\lambda$ such that
$Aut_A(M)\sub G.$ We show that if $M^*$ is a saturated model of an
unsuperstable theory of
cardinality $> Th(M),$ then $M^*$ has the small index property.
\end{abstract}

\section{Introduction}

Throughout the paper we work in ${\frak C}^{eq},$ and we assume that
$M^*$ is a saturated model of
$T$ of cardinality $\lambda.$ We denote the set of
automorphisms of $M^*$ by $Aut(M^*)$ and the set of automorphisms of
$M^*$ fixing $A$ pointwise by $Aut_A(M^*).$  $M^*$ is said to have the small index
property if whenever $G$ is a subgroup of $Aut(M^*)$ with index not larger than
$\lambda$ then for some $A\subset M^*$ with $\big|A\big|<\lambda,$
$Aut_A(M^*)\subseteq G.$ The main theorem of this paper is 
the following result of Shelah: If $M^*$ is a saturated model of
cardinality $\lambda>\big|T\big|$ and there is a tree of height some
uncountable regular cardinal $\kappa\geq\kappa_r(T)$ with
$\mu>\lambda$ 
many branches but at most $\lambda$ nodes, then $M^*$ has the
small index property, in fact 
$$[Aut(M^*):G]\geq \mu$$
for any subgroup $G$ of $Aut(M^*)$ such that for no $A\sub M^*$ with
$|A|<\lambda$ is $Aut_A(M^*)\subseteq G.$ By a result of Shelah on
cardinal arithmetic this implies that if 
$Aut(M^*)$ does not have the small index property, then for some
strong limit $\mu$ such that $cf\,\mu=\aleph_0,$ 
$$\mu<\lambda<2^{\mu}$$
So in particular, if $T$ is unsuperstable,
$M^*$ has the small index property.

\vspace{.1in}

\noindent In the paper ``Uncountable Saturated Structures have the Small Index
Property'' by Lascar and Shelah, the following result was obtained:

\begin{theorem}
Let $M^*$ be a saturated model of cardinality $\lambda$ with
$\lambda>\big|T\big|$ and $\lambda^{<\lambda}=\lambda.$ Then
if $G$ is a subgroup of $Aut(M^*)$ such that 
for no $A\subseteq M^*$
with $\big|A\big|<\lambda$ is $Aut_A(M^*)\subseteq G$ then
$[Aut(M^*):G]$ = $\lambda^{\lambda}.$
\end{theorem}
\proof See [L Sh].

\begin{coro}
Let $M^*$ be a saturated model of cardinality $\lambda$ with
$\lambda>\big|T\big|$ and $\lambda^{<\lambda}=\lambda.$ Then
$M^*$ has the small index property.
\end{coro}

\begin{theorem}
$T$ has a saturated model of cardinality $\lambda$ iff
$\lambda=\lambda^{<\lambda}+D(T)$ or $T$ is stable in $\lambda.$
\end{theorem}
\proof See [Sh c] chp. VIII.

\vspace{.1in}

\noindent So we can assume in the rest of this paper
that $T$ is stable in $\lambda.$ 

\begin{theorem}
$T$ is stable in $\mu$ iff $\mu=\mu_0+\mu^{<\kappa(T)}$ where $\mu_0$
is the first cardinal in which $T$ is stable. 
\end{theorem}
\proof See [Sh c] chp. III.

\vspace{.1in}

\noindent Since $T$ is stable in $\lambda,$ we must have
$\lambda=\lambda^{<\kappa(T)},$ so $cf\,\lambda\geq \kappa(T).$ Since
the first cardinal $\kappa,$ such that $\lambda^{\kappa}>\lambda$ is
regular, we also know that $cf\,\lambda\geq\kappa_r(T).$

\begin{defi}
Let $Tr$ be a tree. If $\eta,\nu\in Tr,$ then $\gamma[\eta,\nu]=$ the least $\gamma$ such
that $\eta(\gamma)\neq \nu(\gamma)$ or else it is
$min(height(\eta),height(\nu)).$
\end{defi}

\begin{notation}
Let $Tr$ be a tree. If $h\in Aut(M^*)$ and $\alpha< height(Tr),\ \eta,\,\nu\in Tr,$ then
$$h^{\eta(\alpha)\,<\,\nu(\alpha)}=h$$
if $\eta(\alpha)<\nu(\alpha)$ and $id_{M^*}$ otherwise.
\end{notation}

\begin{lemma}\label{indep}
Let $\Big\{C_i\mid i\in I\Big\}$ be independent over $A$ and let
$\Big\{D_i\mid i\in I\Big\}$ be independent over $B.$ Suppose that for
each $i\in I$, $tp(C_i/A)$ is stationary. Let $f$ be an elementary map
from $A$ onto $B,$ and let for each $i\in I,$ $f_i$ be an elementary
map extending $f$ which sends $C_i$ onto $D_i.$ Then 
$$\bigcup\limits_{i\in I}f_i$$
is an elementary map from $\bigcup\limits_{i\in I}C_i$ onto $\bigcup\limits_{i\in I}D_i.$
\end{lemma}
\proof Left to the reader.

\begin{lemma}\label{3steps2}
Let $\big|T\big|<\lambda.$ Let $Tr$ be a tree of height $\omega$ with
$\kappa_n$ nodes of 
height $n$ for some $\kappa_n<\lambda.$  Let $n<\omega$ and let
$\langle M_{i}\mid
i\leq n\rangle$ be an increasing chain of models. Let
$M_{n}\subseteq N_0\subseteq N_1\subseteq M^*$ with
$\big|N_1\big|<\lambda.$
Suppose $\langle h_{i}\mid i\leq n\rangle$ are automorphisms of $M^*$ such that

\begin{enumerate}
\item $h_{i}=id_{M_{i}}$
\item $h_{i}[N_j]=N_j$ for $j\leq 1$
\item $h_{i}[M_{k}]=M_{k}$ for $k\leq n$
\end{enumerate}

\noindent For each $\nu\in Tr\restriction level(n+1)$ let
$m_{\nu},l_{\nu}$ be automorphisms of $N_0.$ Let $\eta\in
Tr\restriction level(n+1).$ Suppose $g_{\eta}\in
Aut(N_0)$ such that for all $\nu\in Tr\restriction level(n+1),$
$$g_{\eta}m_{\eta}(m_{\nu})^{-1}(g_{\eta})^{-1}=l_{\eta}(l_{\nu})^{-1}h_{\gamma[\eta,\nu]}^{\eta(\gamma[\eta,\nu])\,<\,\nu(\gamma[\eta,\nu])}$$
Let $m_{\nu}^+,l_{\nu}^+$ be extensions of $m_{\nu}$ and $l_{\nu}$ to
automorphisms of $N_1$ for
all $\nu\in Tr\restriction level(n+1).$ Then there exists a model
$N_2\subseteq  M^*$ containing $N_1$ such that $|N_2|\leq |N_1|+|T|+\kappa_{n+1}$
and $h_i[N_2]=N_2$ for $i\leq n$ and a
$g'_{\eta}\in Aut(N_2)$ extending $g_{\eta}$ and for
all $\nu\in Tr\restriction level(\alpha+1)$ automorphisms of $N_2,$
$m_{\nu}'$ and $l_{\nu}'$ extending $m_{\nu}^+$ and $l_{\nu}^+$
respectively such that
$$g'_{\eta}m'_{\eta}(m'_{\nu})^{-1}(g'_{\eta})^{-1}=l'_{\eta}(l'_{\nu})^{-1}h_{\gamma[\eta,\nu]}^{\eta(\gamma[\eta,\nu])\,<\,\nu(\gamma[\eta,\nu])}$$
\end{lemma}
\proof Let $g^+_{\eta}$ be a map with domain $N_1$ such that
$g^+_{\eta}(N_1)\,\sindep{N_0}\,N_1,$ $g^+_{\eta}(N_1)\subseteq M^*$  and
$g^+_{\eta}$ extends $g_{\eta}.$ Let $g^{++}_{\eta}$ be a
map extending $g_{\eta}$ such that the domain of $(g^{++}_{\eta})^{-1}$ is $N_1,$ 
$(g^{++}_{\eta})^{-1}(N_1)\subseteq M^*$ and
$(g^{++}_{\eta})^{-1}(N_1)\,\sindep{N_0}\,N_1.$ So $g^+_{\eta}\cup g^{++}_{\eta}$ is an
elementary map. Let $l_{\eta}''$ and $m_{\eta}''$ be an extensions of
$l_{\eta}^+$ and $m_{\eta}^+$ to an
automorphisms of $M^*.$ Let
$$m_{\nu}^{++}=
(g_{\eta}^{++})^{-1}(h_{\eta,\nu})^{-1}l_{\nu}(l_{\eta})^{-1}g_{\eta}^{++}m_{\eta}''\restriction
\ (m_{\eta}')^{-1}[(g^{++})^{-1}[N_1]]$$
Note that $m_{\nu}^{+}\cup m_{\nu}^{++}$ is an elementary map. 
Let
$$l_{\nu}^{++}=
(l_{\eta}'')^{-1}g_{\eta}^+m_{\eta}^+(m_{\nu}^+)^{-1}(g_{\eta}^+)^{-1}
(h_{\eta,\nu})^{-1}\restriction
\ h_{\eta,\nu}[g_{\eta}^+[N_1]]$$
Note that $l_{\nu}^{+}\cup l_{\nu}^{++}$ is an elementary map. 
Let $g_{\eta}'',\ m_{\nu}'',\ l_{\nu}''$ be elementary extensions to
$M^*$ of $g^+_{\eta}\cup g^{++}_{\eta},$ $m_{\nu}^{+}\cup
m_{\nu}^{++},$ and $l_{\nu}^{+}\cup l_{\nu}^{++}.$
Let $N_2$ be a model of size $|N_1|+|T|+\kappa_{n+1}$ containing $N_1$ such that
$N_2$ is closed under $m_{\eta}'',\ g_{\eta}'',\ l_{\eta}''$ all the 
$h_{\eta,\nu}$ and $m_{\nu}'',\ l_{\nu}''.$ Let $m_{\nu}',\ l_{\nu}',\
g_{\eta}',\ h_{\eta,\nu}',\ m_{\eta}',\ l_{\eta}'$ be the restrictions to
$N_2$ of the $m_{\nu}'',\ l_{\nu}'',\ g_{\eta}'',\ h_{\eta,\nu}'',\
m_{\eta}'',\ l_{\eta}''.$

\begin{theorem}\label{mine}
If $\lambda>\big|T\big|,$ $cf\,\lambda=\omega,$ $M^*$ is a
saturated model of cardinality $\lambda$ and
if $G$ is a subgroup of $Aut(M^*)$ such that 
for no $A\subseteq M^*$
with $\big|A\big|<\lambda$ is $Aut_A(M^*)\subseteq G$ then
$[Aut(M^*):G]$ = $\lambda^{\omega}.$
\end{theorem}
\proof Suppose not. Let $\big\{\kappa_i\mid i<\omega\big\}$ be an increasing
sequence of cardinals each greater than $\big|T\big|$ with $sup=
\lambda.$ Let $Tr=\big\{\eta\in\, ^{<\omega}\!\lambda\mid
\eta(i)<\kappa_i\big\}.$ Let $M^*=\bigcup\limits_{i<\omega}B_i$ with
$\big|B_i\big|\leq \kappa_i.$ 
By induction  on $n< \omega$ for every $\eta\in
Tr\restriction level\,n$ we define models $N_n\subset M^*$ and $h_n\in
Aut_{N_n}(M^*)-G$ such that
$B_n\sub N_n$ and $\big|N_n\big|\leq\kappa_n,$ and automorphisms
$g_{\eta},m_{\eta},l_{\eta}$ of $N_{n}$ such that
if $\rho\neq\nu$ then $l_{\rho}\neq l_{\nu}$ and
$$g_{\rho}m_{\rho}(m_{\nu})^{-1}(g_{\rho})^{-1}=l_{\rho}(l_{\nu})^{-1}
h_{\gamma[\rho,\nu]}^{\rho(\gamma[\rho,\nu])\,<\,\nu(\gamma[\rho,\nu])}$$
Suppose we have defined the $g_{\eta},m_{\eta},l_{\eta}$ for
$height(\eta)\leq m,$ and $N_j$ for $j\leq m.$  If $n=m+1,$ for each
$i<\kappa_n$ we define models $N_{n,i}$ such that $B_n\sub N_{n,i},\
N_m\sub N_{n,i},\  \langle N_{n,i}\mid i<\kappa_n\rangle$ is increasing
continuous, and for some $\eta_i\in Tr\restriction
level\,n,\ g_{\eta_i}\in Aut(N_{n,i})$ such that
for each $\eta\in Tr\restriction level\,n,\ \eta=\eta_i$ cofinally
many times in $\kappa_n,$ and for every $\nu\in Tr\restriction
level\,n,$ $m_{\nu}^i\neq l_{\nu}^i\in Aut(N_{n,i})$ such that
$$g_{\eta_i}m_{\eta_i}^i(m_{\nu}^i)^{-1}(g_{\eta_i})^{-1}=l_{\eta_i}^i(l_{\nu}^i)^{-1} 
h_{\gamma[\eta_i,\nu]}^{\eta_i(\gamma[\eta_i,\nu])\,<\,\nu(\gamma[\eta_i,\nu])}$$
The $g_{\eta_i},\ m_{\nu}^i,\ l_{\nu}^i$ are easily defined by
induction on $i<\kappa_n$ using lemma \ref{3steps2} so that if
$i_1<i_2$ then $m_{\nu}^{i_1}\subseteq m_{\nu}^{i_2},\
l_{\nu}^{i_1}\subseteq l_{\nu}^{i_2},$ and if $\eta_{i_1}=\eta_{i_2}$
then $g_{\eta_{i_1}}\subseteq g_{\eta_{i_2}}.$ 
Then if we let 
$g_{\eta}=\bigcup\big\{g_{\eta_i}\mid\eta_i=\eta\big\},\
m_{\eta}=\bigcup\limits_{i<\kappa_n}m_{\eta}^i,$ 
$l_{\eta}=\bigcup\limits_{i<\kappa_n}l_{\eta}^i,$ 
$N_n=\bigcup\limits_{i<\kappa_n}N_{n,i}$ and $h_n\in Aut_{N_n}(M^*)-G$  we have finished. Let
$Br$ be the set of branches of $Tr$ of height $\omega.$ For
$\rho\in Br$ let $g_{\rho}=\bigcup\big\{g_{\eta}\mid \eta<\rho\big\},\
m_{\rho}=\bigcup\big\{m_{\eta}\mid \eta<\rho\big\},$ and $ 
l_{\rho}=\bigcup\big\{l_{\eta}\mid \eta<\rho\big\}.$ If $\rho\neq\nu,$
$g_{\rho}\neq g_{\nu}$ since without loss of generality
$\rho(\gamma[\rho,\nu])<\nu(\gamma[\rho,\nu])$ and
$$g_{\rho}m_{\rho}(m_{\nu})^{-1}(g_{\rho})^{-1}=l_{\rho}(l_{\nu})^{-1}
h_{\gamma[\rho,\nu]}^{\rho(\gamma[\rho,\nu])\,<\,\nu(\gamma[\rho,\nu])}$$
and
$$g_{\nu}m_{\nu}(m_{\rho})^{-1}(g_{\nu})^{-1}=l_{\nu}(l_{\rho})^{-1}$$
implies
$$g_{\rho}(g_{\nu})^{-1}l_{\rho}(l_{\nu})^{-1}g_{\nu}(g_{\rho})^{-1}=
l_{\rho}(l_{\nu})^{-1}
h_{\gamma[\rho,\nu]}^{\rho(\gamma[\rho,\nu])\,<\,\nu(\gamma[\rho,\nu])}$$
So if $g_{\rho}=g_{\nu}$ this would imply
$h_{\gamma[\rho,\nu]}^{\rho(\gamma[\rho,\nu])\,<\,\nu(\gamma[\rho,\nu])}=id_{M^*}$
a contradiction. If
$$[Aut(M^*):G]<\lambda^{\omega}$$
then for some $\rho,\nu\in Br$ we must have $l_{\rho}(l_{\nu})^{-1}\in
G$ and $g_{\rho}(g_{\nu})^{-1}\in G,$ but then we get a contradiction
as 
$g_{\rho}(g_{\nu})^{-1}l_{\rho}(l_{\nu})^{-1}g_{\nu}(g_{\rho})^{-1}\in
G$ and $l_{\rho}(l_{\nu})^{-1}\in G,$ but
$h_{\gamma[\rho,\nu]}^{\rho(\gamma[\rho,\nu])\,<\,\nu(\gamma[\rho,\nu])}\not\in
G.$

\begin{coro}
If $\lambda>\big|T\big|,$ $cf\,\lambda=\omega$ and $M^*$ is a
saturated model of cardinality $\lambda$ then $M^*$ has the small
index property.
\end{coro}

\vspace{.1in}

\noindent So we will assume in the remainder of the paper that in
addition to $T$ being stable, $cf\,\lambda\geq\kappa_r(T)+\aleph_1$
and 
$T,M^*,$ and $\lambda$ are constant.

\section{Constructing $M^*$ as a chain from $K_{\delta}$}

\begin{defi}
Let $\delta<\lambda^+,\ cf\,\delta\geq \kappa_r(T).$ 
$$K^s_{\delta}=\Big\{\bar{N}\mid \bar{N}=\langle N_i\mid 
i\leq
\delta\rangle,\  N_i\hbox{ is increasing continuous, } |N_i|=\lambda,$$
$$N_0\hbox{ is
saturated, }N_{\delta}=M^*,\hbox{ and }\big(N_{i+1},c\big)_{c\in
N_i}\hbox{ is saturated}\Big\}$$
 For $\mu>\aleph_0,$ 
$$K^{\mu}_{\delta}=\Big\{\bar{A}\mid \bar{A}=\langle A_i\mid
 i\leq
\delta\rangle,$$
$$A_i\hbox{ is increasing continuous, }
\big|A_{\delta}\big|<\mu,\ acl\,A_i=A_i\Big\}$$
If $\bar{A}\in K^{\lambda^+}_{\delta},$ then $f\in Aut(\bar{A})$ if $f$
is an elementary permutation of $A_{\delta}$ and  if $i\leq
\delta,$ then $f\restriction A_i$ is a permutation of $A_i.$ 
\end{defi}

\begin{defi}
Let $\bar{A}^0,\bar{A}^1\in K^{\mu}_{\delta}.$ Then $\bar{A}^0\leq
\bar{A}^1$  iff   $\bigwedge\limits_{i\leq\delta}A_i^0\subseteq A_i^1$ and
$i<j\leq  \delta\ \ \Rightarrow\ \ A_i^1\,\sindep{A_i^0}\,A_j^0.$ 
\end{defi}

\begin{lemma} 
\begin{enumerate}
\item $\big(K^{\mu}_{\delta},\leq\big)$ is a partial order
\item Let $\bar{A}^\zeta\in K^{\mu}_{\delta}$ for $\zeta<\zeta(*)$
and let $\xi<\zeta\ \Rightarrow\ \bar{A}^{\xi}\leq \bar{A}^{\zeta}.$ If we let
$A_i=\bigcup_{\zeta<\zeta(*)} A_i^{\zeta},$ and
$\Big|\bigcup_{\zeta<\zeta(*)}A_i^{\delta}\Big|<\mu,$ then 
$$\bar{A}=\langle A_i\mid i\leq \delta\rangle\in K^{\mu}_{\delta}$$
and for every $\zeta<\zeta(*),\  \bar{A}_{\zeta}\leq \bar{A}.$
\item If $\bar{A}^{\zeta}\leq \bar{A}^*$ for $\zeta<\zeta(*),$ and
$\bar{A}$ is as above, then
$\bar{A}\leq \bar{A}^*$
\end{enumerate}
\end{lemma}
\proof 
\begin{enumerate}
\item By the transitivity of nonforking.
\item By the finite character of forking.
\item By the finite character of forking.
\end{enumerate}

\begin{defi}
Let $A\subseteq M,$ with
$\big|A\big|<\kappa_r(T)$
and let $p\in
S(acl\,A).$ Then $dim(p,M)=$ the
minimal cardinality of an maximal independent set of realizations of
$p$ inside $M.$ If $M$ is $\kappa_r^{\epsilon}(T)$-saturated
($\kappa_r^{\epsilon}$-saturated means $\aleph_{\epsilon}$-saturated
if $\kappa_r(T)=\aleph_0$ and $\kappa_r(T)$ saturated otherwise) then by
[Sh c] III 3.9. $dim(p,M)=$the cardinality of any maximal independent
set of realizations of $p$ inside $M.$ 
\end{defi}

\begin{lemma}\label{dimsat}
Let $\big|M\big|=\lambda$ and assume that $M$ is 
$\kappa_r^{\epsilon}(T)$-saturated. Then $M$ is saturated 
if and only if for every $A\subseteq M,$ with
$\big|A\big|<\kappa_r(T)$
and  $p\in
S(acl\,A),$ $dim(p,M)=\lambda.$
\end{lemma}
\proof See [Sh c] III 3.10.

\begin{lemma}\label{insat}
Let $\langle \bar{A}^{\alpha}\mid \alpha<\lambda\rangle$ be an
increasing continuous sequence of elements of
$K^{\lambda^+}_{\delta}$ such that
$\forall\,\gamma<\delta,\ \forall\,A\subseteq
\bigcup\limits_{\alpha<\lambda}A^{\alpha}_{\gamma}$ if
$\big|A\big|<\kappa_r(T)$ and $p\in S(acl\,A)$ then for $\lambda$ many
$\alpha<\lambda,$

\begin{enumerate}
\item $A_{\zeta}^{\alpha}=A^{\alpha+1}_{\zeta}\
\forall\,\zeta\leq\gamma$
\item There exists $a\in A_{\gamma+1}^{\alpha+1}$ such that the type of
$a/A_{\gamma+1}^{\alpha}$ is the stationarization of $p$
\end{enumerate}

\noindent then 
$$\langle N_{\gamma}\mid \gamma<\delta\rangle\in
K^s_{\delta}$$
where $N_{\gamma}=\bigcup\limits_{\alpha<\lambda}A_{\gamma}^{\alpha}.$
\end{lemma}
\proof It is enough to show $\forall\,\gamma<\delta$ that $\big(N_{\gamma+1},\ c\big)_{c\,\in\,N_{\gamma}}$ is
saturated. For this by lemma \ref{dimsat} it is enough to show
$\forall\,A\subseteq N_{\gamma+1}$ such that $\big|A\big|<\kappa_r(T)$
and for every type $p\in S(acl\,A\cup N_{\gamma}),$
$$dim (p, N_{\gamma+1})=\lambda$$
By the
assumption of the lemma, there exists $\big\{a_i\mid i<\lambda\big\}$
realizations of $p\restriction acl\,A$ and
$\langle A_{\gamma+1}^{\alpha_i}\mid i<\lambda\rangle$ 
such that for each $i<\lambda,\ a_i\in A_{\gamma+1}^{\alpha_i+1},
\ A_{\gamma}^{\alpha_i+1}=A_{\gamma}^{\alpha_i},$ and
$$a_i\,\indep{A}\,A_{\gamma+1}^{\alpha_i}\ \ \hbox{ and }\ \
a_iA_{\gamma+1}^{\alpha_i}\,\indep{A_{\gamma}^{\alpha_i}}\,N_{\gamma}$$
which implies
$$a_i\,\indep{A_{\gamma+1}^{\alpha_i}}\,N_{\gamma}\ \ \hbox{ and }\ \
a_i\,\indep{A}\,N_{\gamma}$$
Since $cf\,\lambda\geq\kappa_r(T)$ without loss of generality $A\subseteq
A_{\gamma+1}^{\alpha_0}.$ We must show the $\langle a_i\mid i<\lambda\rangle$ are independent
over $N_{\gamma}\cup A.$ By induction on $i<\lambda,$ we show that
$$\langle a_j\mid j\leq i\rangle$$
are independent over $A\cup \big\{A_{\gamma}^{\alpha_j}\mid
j\leq i\big\}.$ This is enough as
$$\big\{a_j\mid j\leq i\big\}\,\indep{A\cup
\big\{A_{\gamma}^{\alpha_j}\mid j\leq i\big\}}\,N_{\gamma}$$
Since $\langle a_j\mid j<i\rangle$ are independent over $A\cup
\big\{A_{\gamma}^{\alpha_i}\mid j<i\big\},$ and 
$$\big\{a_j\mid j<i\big\}\,\indep{A\cup
\big\{A_{\gamma}^{\alpha_j}\mid j<i\big\}}\,A_{\gamma}^{\alpha_i}$$
$\langle a_j\mid j<i\rangle$ are independent over $A\cup
A_{\gamma}^{\alpha_i}.$
Since $a_i\,\sindep{A\cup A_{\gamma}^{\alpha_i}}\,A_{\gamma+1}^{\alpha_i}$
we have
$$a_i\,\indep{A\cup A_{\gamma}^{\alpha_i}}\,\big\{a_j\mid j<i\big\}$$

\begin{lemma}
Let $\langle \bar{N}^{\alpha}\mid\alpha<\delta\rangle$ be an
increasing continuous sequence of elements of $K^{\mu^+}_{\delta}$ such that
$\bigcup\limits_{\alpha<\delta}N_{\delta}^{\alpha}=M^*$ and for every $\gamma<\delta,$ and $\alpha<\delta,$ 
$$\big(N_{\gamma+1}^{\alpha+1},c\big)_{c\in N_{\gamma+1}^{\alpha}\cup
N_{\gamma}^{\alpha+1}}$$
and
$$(N_0^{\alpha+1},c)_{c\in N_0^{\alpha}}$$
are saturated of cardinality $\lambda.$ Then 
$$\langle N_{\gamma}\mid \alpha<\delta\rangle\in K^s_{\delta}$$
where $N_{\gamma}=\bigcup\limits_{\alpha<\delta}N_{\gamma}^{\alpha}.$
\end{lemma}
\proof Similar to the proof of the previous lemma.

\begin{lemma}\label{build4}
Let $cf\,\delta\geq \kappa_r(T)+\aleph_1.$ 
Let $\bar{M}\in K^s_{\delta}.$ Let
$A_{\delta}\subseteq M^*$ such that
$\big|A_{\delta}\big|<\lambda$ and
$A_{\delta}=\bigcup\limits_{i<\delta}A_i$
where $\langle A_i\mid i<\delta\rangle$ is an increasing
continuous chain. Suppose $\forall\beta<\delta,$ and
$\forall i<\delta,$
$$M_{\beta}\,\indep{A_i\cap M_{\beta}}\,A_i$$
Let $a\subseteq M_{\beta^*}$ such
that $\big|a\big|<\kappa_r(T).$ Then there exists a continuous
increasing sequence  $\langle A_i'\mid i<\delta\rangle$ and a set $B$ such
that $\big|B\big|<\kappa_r(T),$ $A_i\subset A_i',\ a\subset\bigcup A_i'=A_{\delta}',\
\big|A_{\delta}'\big|<\lambda,$ for some non-limit $i^*<\delta,$ $A_i'=A_i$ if
$i<i^*,$ and $A_i'=A_i\cup B$ if $i^*\leq i$ and $\forall\ i,\beta<\delta,$
$$M_{\beta}\,\indep{A_i'\cap M_{\beta}}\,A_i'$$
and $\forall\ i,\beta<\delta,$
$$M_{\beta}\cup (M_{\beta+1}\cap A_{\delta})\,\indep{M_{\beta}\cup (M_{\beta+1}\cap
A_i)}\,A_i'\cap M_{\beta+1}$$
and
$$A_{\delta}\,\indep{A_i}\,A_i'$$
\end{lemma}
\proof First by induction on $n\in\omega,$ we define $\langle B_n\mid
n<\omega\rangle$ such that $B_0=a,$ $\big|B_n\big|<\kappa_r(T)$ and
$\forall\ i<\delta,\ \forall\beta<\delta,$
$$B_n\,\indep{(M_{\beta}\cap(A_i\cup B_{n+1}))\cup A_i}\,M_{\beta}\cup A_i$$
So suppose $B_n$ has been defined. By induction on $m<\omega$ we define
subsets $C_1$ and $C_2$ of $\delta$ such that $0\in C_i,$ $|C_i|\leq
\kappa_r(T)$ and such that if $(a_1,b_1),(a_2,b_1),(a_1,b_2),(a_2,b_2)$
are four neighboring points in $C_1\times C_2$ with $a_1<a_2$ and $b_1<b_2,$
then for all $i,j$ such that $a_1\leq i<a_2$ and $b_1\leq j<b_2$
$$B_n\,\indep{M_{a_1}\cup A_{b_1}}\,M_{a_1+i}\cup A_{b_1+j}$$
So it is enough to find $\big|B_{n+1}\big|<\kappa_r(T)$ such that for
every $(a,b)\in C_1\times C_2,$
$$B_n\,\indep{(M_a\cap (A_b\cup B_{n+1}))\cup A_b}\,M_a\cup A_b$$
As $\big|C_1\times C_2\big|<\kappa_r(T)$ this is possible. 
Let $B=\bigcup\limits_{n\in\omega}B_n.$ (If $\kappa_r(T)=\aleph_0$
then without loss of generality we can define the $B_n$ such that for
some $k<\omega,\ \bigcup\limits_{n\in \omega}B_n=\bigcup\limits_{n\in
k}B_n.)$ It is enough to prove the
following statement. 

\vspace{.1in}

\noindent {\em There exists a non-limit $i^*<\delta$ such that if
$A_i'=A_i$ for 
$i<i^*,$ and $A_i'=A_i\cup B$ for $i\geq i^*$ then the conditions of
the theorem hold.}

\noindent \proof $\forall\beta<\delta,\ \forall i<\delta,$ if
$A_i'=A_i\cup B,$ then since
$$B\,\indep{(M_{\beta}\cap(A_i\cup B))\cup A_i}\,M_{\beta}\cup A_i$$
we have
$$A_i'\,\indep{A_i'\cap M_{\beta}}\,M_{\beta}$$
Let $i^{**}<\delta$ such that for all $i\geq i^{**},$
$$A_{\delta}\,\indep{A_i}\,A_i'$$
It is enough to find $i^{**}\leq i^*<\delta$ such that $\forall\beta<\delta,$
$$B\,\indep{M_{\beta}\cup (M_{\beta+1}\cap A_{i^*})}\,M_{\beta}\cup (M_{\beta+1}\cap A_{\delta})$$
Let $\langle \beta_{\alpha}\mid \alpha\in\gamma\rangle$ where
$\gamma<\kappa_r(T)$ be the set of all places such that
$$B\,\noindep{M_{\beta_{\alpha}-1}\cup (M_{\beta_{\alpha}}\cap
A_{\delta})}\,M_{\beta_{\alpha}}\cup (M_{\beta_{\alpha}+1}\cap A_{\delta})$$
For each $\beta\in \langle \beta_{\alpha}\mid \alpha\in\gamma\rangle$
let $i_{\alpha}$ be such that
$$B\,\indep{M_{\beta_{\alpha}}\cup (M_{\beta_{\alpha}+1}\cap
A_{i_{\alpha}})}\,M_{\beta-1}\cup (M_{\beta_{\alpha}+1}\cap A_{\delta})$$
Let $i_{\gamma}$ be such that
$$B\,\indep{M_{0}\cup (M_{1}\cap
A_{i_{\gamma}})}\,M_{0}\cup (M_{1}\cap A_{\delta})$$
Let $i^*=sup\big\{i_{\alpha}\mid \alpha\in\gamma+1\big\}+1+i^{**}.$
As $\big|B\big|<\kappa_r(T)$ and $cf\,\delta\geq\kappa_r(T),$
$i^*<\delta,$ so there is
no problem.

\begin{lemma}\label{build6}
Let $\bar{M}\in K^s_{\delta}.$ Let $A\subseteq M^*$ such that $\big|A\big|<\lambda$ and
$A=\bigcup\limits_{i<\delta}A_i$ where $\langle
A_i\mid i<\delta\rangle$ is increasing continuous,
each $A_i$ is algebraically closed and $\forall\,i <\delta,\ \forall\,\beta<\delta,$
$$M_{\beta}\,\indep{M_{\beta}\cap A_i}\,A_i$$
Let $i^*$ be a successor $<\delta,\ \beta^*<\delta,
\ \beta^*$ a successor, and let $p\in S(A_{i^*}\cap M_{\beta^*}).$ (Or
even a $<\lambda$ type over $A_i\cap M_{\beta^*}.$) Let $p'\in
S\big((A_{i^*}\cap M_{\beta^*})\cup M_{\beta^*-1}\big)$ such that $p'$ does
not fork over $p.$ Then there exists an $a\in M_{\beta^*}$ such that
$a$ realizes $p',$
$$A\,\indep{M_{\beta^*}\cap
A_{i^*}}\,a$$
and if $A_i'$ = $A_i\cup\big\{a\big\}$ for $i\geq
i^*$ and $A_i'=A_i$ for $i<i^*,$ then $\forall\,\beta<\delta,\ \forall\,i<\delta,$
$$M_{\beta}\,\indep{M_{\beta}\cap A_i'}\,A_i'$$
\end{lemma}
\proof Let $B\subseteq M_{\beta^*}$ such that $\big|B\big|<\lambda,\ A_i^*\cap
M_{\beta^*}\subseteq B,$ and 
$$M_{\beta^*}\,\indep{M_{\beta^*-1}B}\,A$$
Let $a\in M_{\beta^*}$ such that $a$ realizes  $p$ and
$$a\,\indep{A_i^*\cap M_{\beta^*}}\,B\cup M_{\beta^*-1}$$
Since
$$M_{\beta^*}\,\indep{M_{\beta^*-1}\cup B}\,A$$ 
we have
$$a\,\indep{M_{\beta^*-1}\cup B}\,A$$  
which implies 
$$a\,\indep{A_i^*\cap
M_{\beta^*}}\,M_{\beta^*-1}\cup A$$ 
Since for all $i\geq i^*,$ 
$$a\,\indep{A_i}\,M_{\beta^*-1}\cup A$$
we have for all $\gamma<\beta^*,$
$$a\,\indep{A_i}\,M_{\gamma}\cup A$$
which implies
$$a\cup A_i\,\indep{A_i\cap M_{\gamma}}\,M_{\gamma}$$
Since $a\subseteq  M_{\beta^*}$ we also have $\forall\gamma\geq
\beta^*,$
$$a\cup A_i\,\indep{(a\cup A_i)\cap M_{\gamma}}\,M_{\gamma}$$

\begin{lemma}\label{build7}
Let $\bar{M}\in K^s_{\delta}.$ Let $A\subseteq M^*$ such that $\big|A\big|<\lambda$ and
$A=\bigcup\limits_{i<\delta}A_i$ where $\langle
A_i\mid i<\delta\rangle$ is increasing continuous,
each $A_i$ is algebraically closed and $\forall\,i <\delta,\  \forall\,\beta<\delta,$
$$M_{\beta}\,\indep{M_{\beta}\cap A_i}\,A_i$$
Let $i^*<\delta,\  \beta^*<\delta,
\ \beta^*,\ i^*$ successors, and let $p\in S(A_i\cap M_{\beta}).$ Let $p'\in
S\big((A_i\cap M_{\beta^*})\cup M_{\beta^*-1}\big)$ such that $p'$ does
not fork over $p.$ Let $f\in Aut(A)$ such that
$\forall\,i<\delta,\ f[A_i]=A_i.$ Then there exists 
$\big\{a_i\mid i\in \Bbb{Z}\big\}\subseteq M^*$ and
an extension $f'$ of $f$ with domain $A\cup  \big\{a_i\mid i\in
\Bbb{Z}\big\}$ such that $a_0$ realizes $p',\ a_0\in M_{\beta^*},$ and $\forall\,i\in
\Bbb{Z}\ f'(a_i)=a_{i+1}$ and if $A_i'=A_i\cup \big\{a_i\mid i\in
\Bbb{Z}\big\}$ for $i\geq i^*$ and $A_i'=A_i$ for $i<i^*,$ then for all
$\beta<\delta,$ 
$$M_{\beta}\,\indep{M_{\beta}\cap A_i'}\,A_i'$$
$$A_{\delta}\,\indep{A_i}\,A_i'$$
and
$$M_{\beta-1}\cup (M_{\beta}\cap A)\,\indep{M_{\beta-1}\cup (M_{\beta}\cap A_i)}\,M_{\beta}\cap A_i'$$  
\end{lemma}
\proof We define $\big\{a_i\mid i\in -n,\ldots,0,\ldots,n\big\}$ by
induction on $n$ such that if $A_i'=acl(A_i\cup \big\{a_i\mid i\in
-n,\ldots,0,\ldots,n\big\})$ if $i\geq i^*$ and $A_i'=A_i$ if $i<i^*,$
then
$\forall i<\delta,\ \forall\beta<\delta,$
$$M_{\beta}\,\indep{M_{\beta}\cap A_i'}\,A_i'$$
$$A_{\delta}\,\indep{A_i}\,A_i'$$
and
$$M_{\beta-1}\cup (M_{\beta}\cap A)\,\indep{M_{\beta-1}\cup (M_{\beta}\cap A)}\,M_{\beta}\cap A_i'$$
and $f_n=f\cup \big\{(a_i,a_{i+1})\mid -n\leq i<n\big\}$ is an elementary
map. 
In addition we define a sequence of successor ordinals $\langle
\beta_i\mid i\in{\Bbb Z}\rangle$ such that $\beta_i<\beta_j$ if
$|i|<|j|,$ and $\beta_n<\beta_{-n}$ such that
$$a_{n+1}\,\indep{M_{\beta_{n+1}}\cap
A_{i^*}}\,M_{\beta_{n+1}-1}\cup A\cup\{a_{-n}\ldots,a_0,\ldots,a_n\}$$
and
$$a_{-(n+1)}\,\indep{M_{\beta_{-(n+1)}}\cap
A_{i^*}}\,M_{\beta_{-(n+1)}-1}\cup A\cup\{a_{-n},\ldots,a_0,\ldots,a_n,a_{n+1}\}$$
Define $a_0$ as in the previous lemma. Suppose that
$\big\{a_{-n},\ldots,a_0,\ldots,a_n\big\}$ and $\beta_i$ for $-n\leq
i\leq n$ have been defined satisfying the conditions. Let $C= acl\,C$
such that for some $B\subseteq C$ with $\big|B\big|<\kappa_r(T),$ 
$acl B=C,$ $C\subseteq M_{\beta_{-n}}\cap A_{i^*}$ and
$$a_n\,\indep{C}\,A\cup\big\{a_{-n},\ldots,a_0,\ldots,a_{n-1}\big\}$$
Let $\beta_{n+1}>\beta_{-n}$ be a successor such that $f(C)\subseteq
M_{\beta_{n+1}}\cap A_{i^*}.$ Let $a_{n+1}\in M_{\beta_{n+1}}$ realize
$$f_n\Big(tp\big(a_n/A\cup\big\{a_{-n},\ldots,a_0,\ldots,a_{n-1}\big\}\big)\Big)$$
and in addition 
$$a_{n+1}\,\indep{M_{\beta_{n+1}}\cap A_i^*}\,A\cup M_{\beta_{n+1}-1}$$
Similarly for $a_{-(n+1)}.$ Now as in the proof of the previous lemma,
all the conditions of the induction hold.

\begin{lemma}\label{build8}
Let $\delta$ be an ordinal less than $\lambda^+$ such that $cf\,\delta\geq\aleph_1+\kappa_r(T).$ Let $f\in Aut_E(M^*)$ with
$\big|E\big|<\lambda.$ 
Let $\bar{M}\in
K^s_{\delta}.$ Then there exists
$\bar{N}^1,\bar{N}^2\in K^s_{\delta},\ f_1\in
Aut_E(\bar{N}^1),\ f_2\in
Aut_E(\bar{N}^2)$ with $E\subseteq N^1_0,\ E\subseteq N^2_0$ such that

\begin{enumerate}
\item $f=f_2f_1$
\item $\forall\,i,\beta<\delta,
\ \forall\,l\in\{0,1\},$
$$M_{\beta}\,\indep{M_{\beta}\cap N^l_i}\,N^l_i$$
\item $\forall\,i,\beta<\delta,
\ \forall\,l\in\{0,1\},$
$$\big(N^l_{i+1}\cap M_{\beta+1}, c\big)_{c\,\in\,(N_{i+1}^l\cap
M_{\beta})\,\cup\,(N^l_i\cap M_{\beta+1})}$$ is saturated of cardinality
$\lambda$
\item $\big(N_{i+1}^l\cap M_0, c\big)_{c\,\in\,N_{i}^l\,\cap\,M_0}$
is saturated of cardinality $\lambda$
\end{enumerate}

\end{lemma}
\proof Without loss of generality $E=\emptyset.$ By induction on
$\alpha<\lambda$ we build increasing continuous 
sequences $\langle A_i^{\alpha}\mid i\leq \delta\rangle,\ \langle B_i^{\alpha}\mid i\leq \delta\rangle,\ \langle f_1^{\alpha}\mid \alpha<\lambda\rangle,\
\langle f_2^{\alpha}\mid \alpha<\lambda\rangle$ such that

\begin{enumerate}
\item $M^*=\bigcup\limits_{\alpha<\lambda}A^{\alpha}_{\delta}=\bigcup\limits_{\alpha<\lambda}B^{\alpha}_{\delta}$  
\item $N^1_i=\bigcup\limits_{\alpha<\lambda}A_i^{\alpha}\ \ \
N^2_i=\bigcup\limits_{\alpha<\lambda}B_i^{\alpha}$
\item $f_1^{\alpha}\in Aut(A_{\delta}^{\alpha})$
such that $f_1^{\alpha}[A_i^{\alpha}]=A_i^{\alpha}$
\item $f_2^{\alpha}\in Aut(B_{\delta}^{\alpha})$
such that $f^{\alpha}_2[B_i^{\alpha}]=B_i^{\alpha}$
\item $f[A_i^{\alpha}]=A_i^{\alpha},\ \ \ f[B_i^{\alpha}]=B_i^{\alpha}$
\item $\big|A_{\delta}^{\alpha}\big|<\big|\alpha\big|^++\kappa_r(T)+\aleph_1$
\item $\big|B_{\delta}^{\alpha}\big|<\big|\alpha\big|^++\kappa_r(T)+\aleph_1$
\item $A_{\delta}^{\alpha}=B_{\delta}^{\alpha}$
\item $f_{\alpha}^2f_{\alpha}^1=f\restriction A_{\delta}^{\alpha}$
\item $\forall\,\beta<\delta,\
\forall\,i<\delta,\ \forall\,\alpha<\lambda,$
$$M_{\beta}\,\indep{M_{\beta}\cap A_i^{\alpha}}\,A_i^{\alpha}$$
\item $\forall\,\beta<\delta,\
\forall\,i<\delta,\ \forall\,\alpha<\lambda,$
$$M_{\beta}\,\indep{M_{\beta}\cap B_i^{\alpha}}\,B_i^{\alpha}$$
\item $\forall\,i,\beta<\delta,
\ \forall\,l\in\{0,1\},$
$$\big(N^l_{i+1}\cap M_{\beta+1}, c\big)_{c\,\in\,(N_{i+1}^l\cap
M_{\beta})\,\cup\,(N^l_i\cap M_{\beta+1})}$$ is saturated of cardinality
$\lambda$
\item $\big(N_{i+1}^l\cap M_0, c\big)_{c\,\in\,N_{i}^l\,\cap\,M_0}$
is saturated of cardinality $\lambda$
\item $\forall i<\delta,\ \forall\alpha<\lambda,$
$$A_{\delta}^{\alpha}\,\indep{A_i^{\alpha}}\,A_i^{\alpha+1}$$
\item $\forall i<\delta,\ \forall\alpha<\lambda,$
$$B_{\delta}^{\alpha}\,\indep{B_i^{\alpha}}\,B_i^{\alpha+1}$$
\item $\forall\beta<\delta,\ \forall i<\delta, \forall \alpha<\lambda,$
$$M_{\beta}\cup (M_{\beta+1}\cap
A_{\delta}^{\alpha})\,\indep{M_{\beta}\cup (M_{\beta+1}\cap
A_i^{\alpha})}\,M_{\beta+1}\cap A_{i}^{\alpha+1}$$
\item $\forall\beta<\delta,\ \forall i<\delta, \forall \alpha<\lambda,$
$$M_{\beta}\cup (M_{\beta+1}\cap
B_{\delta}^{\alpha})\,\indep{M_{\beta}\cup (M_{\beta+1}\cap
B_i^{\alpha})}\,M_{\beta+1}\cap B_{i}^{\alpha+1}$$
\end{enumerate}

\noindent At limit stages we take unions. Let $\alpha$ be even. Let
$M^*=\langle m_{\alpha}\mid \alpha<\lambda\rangle.$ In the induction
we define $\langle p_{\alpha}\mid \alpha$ is even
and $\alpha<\lambda\rangle$ such that each $p_{\alpha}\in
S((M_{\beta+1}\cap A_{i+1}^{\alpha})\cup M_{\beta})$ for some $i,\beta<\delta$
and  such that
$\forall\,i<\delta,\ \forall\,\beta<\delta,\ \forall\,A\subseteq M^*$ such that
$\big|A\big|<\kappa_r(T),\ \forall\,p\in S(acl\,A)$ there exists $\lambda$
many $p_{\alpha}\in \langle p_{\alpha}\mid\alpha<\lambda\rangle$ such
that $p_{\alpha}\in S((M_{\beta+1}\cap A_{i+1}^{\alpha})\cup M_{\beta}),$ $p_{\alpha}$ is
a nonforking extension of $p,$ $p_{\alpha}$ is realized in
$A_{i+1}^{\alpha+1}\cap M_{\beta+1},$ and $\forall j\leq i,
A_j^{\alpha}=A_j^{\alpha+1}.$  By the proof of lemma \ref{insat} this
insures 12. and 13. holds for $l=1$ when we finish our construction.
So let $i^*,\beta^* <\delta$ such that $p_{\alpha}\in
S((M_{\beta^*+1}\cap A_{{i+1}^*}^{\alpha})\cup M_{\beta^*}).$ By lemma
\ref{build7} we can find 
an extensions $(A_{i}^{\alpha})'$ of $A_{i}^{\alpha}$ with
$(A_i^{\alpha})'=A_i^{\alpha}$ for $i\leq i^*$ and extension $f_1'$
of $f_1$ such that $f_1'[(A_{i}^{\alpha})']=(A_{i}^{\alpha})',\
p_{\alpha}$ is realized in $M_{\beta^*+1}\cap (A_{i^*+1}^{\alpha})'$ and
$\forall\,\beta<\delta,\ \forall\,i<\delta,$ 
$$M_{\beta-1}\cup (M_{\beta}\cap
A^{\alpha}_{\delta})\,\indep{M_{\beta-1}\cup (M_{\beta}\cap A_i^{\alpha})}\,M_{\beta}\cap
(A_i^{\alpha})'$$ 
$$A_{\delta}^{\alpha}\,\indep{A_i^{\alpha}}\,(A_i^{\alpha})'$$
and
$$M_{\beta}\,\indep{M_{\beta}\cap (A_i^{\alpha})'}\,(A_i^{\alpha})'$$
Let $F_1'$ be an extension of $f_1'$ to an automorphism of $M^*.$ 
By iterating $\omega$ times the procedure in the proof of lemma
\ref{build4} we can find $D\subset M^*$ such that
$\big|D\big|<\kappa_r(T)+\omega_1,$ if $m$ is the least element of $\langle
m_{\alpha}\mid \alpha<\lambda\rangle$ then $m\in D,$ $D$ is closed
under $f,f^{-1},F_1',(F_1')^{-1}$ and for some $i^{**},i^{***}<\delta$
if $A_i^{\alpha+1}=(A_i^{\alpha})'\cup D,$ for $i\geq i^{**}$ and
$(A_i^{\alpha})'$ for $i<i^{**}$ and if
$B_i^{\alpha+1}=B_i^{\alpha}\cup D,$ for $i\geq i^{***}$ and 
$B_i^{\alpha}$ for $i<i^{***}$
then
$$M_{\beta}\cup (M_{\beta}\cap
A^{\alpha}_{\delta})\,\indep{M_{\beta}\cup (M_{\beta+1}\cap A_i^{\alpha})}\,M_{\beta+1}\cap
A_i^{\alpha+1}$$ 
$$A_{\delta}^{\alpha}\,\indep{A_i^{\alpha}}\,A_i^{\alpha+1}$$
and
$$M_{\beta}\,\indep{M_{\beta}\cap A_i^{\alpha+1}}\,A_i^{\alpha+1}$$
and
$$M_{\beta}\cup (M_{\beta+1}\cap
B^{\alpha}_{\delta})\,\indep{M_{\beta}\cup (M_{\beta+1}\cap B_i^{\alpha})}\,M_{\beta+1}\cap
B_i^{\alpha+1}$$ 
$$B_{\delta}^{\alpha}\,\indep{B_i^{\alpha}}\,B_i^{\alpha+1}$$
and
$$M_{\beta}\,\indep{M_{\beta}\cap B_i^{\alpha+1}}\,B_i^{\alpha+1}$$
Similarily for $\alpha$ odd.
Let $f_1^{\alpha+1}=F_1\restriction A_i^{\alpha+1}$ and
$f_2^{\alpha+1}=f(f_1^{\alpha+1})^{-1}.$

\section{The proof of the small index property}

\begin{defi}
Let $\delta$ be a limit ordinal and let $\bar{N}\in K_{\delta}^s.$ Then $f\in
Aut^*(\bar{N})$ if and only if $f\in Aut(M^*)$ and for some $n\in\omega,$
$f[N_{\alpha}]=N_{\alpha}$ for every $\alpha$ such that
$n\leq\alpha\leq\delta.$ $Aut_A^*(\bar{N})=\big\{f\in
Aut^*(\bar{N})\mid f\restriction A=id_A\big\}.$
\end{defi}

\begin{defi}
Let $\delta$ be a limit ordinal and let $\bar{N}\in K_{\delta}^s.$ Let $B\subseteq N_0$
as in the above definition. If for every
$f\in Aut(M^*)$
$$(\,f\in Aut^{*}(\bar{N})\ \wedge\ f\restriction B=id_{B}\,)\
\Rightarrow \ f\in G$$
then we define
$$E=\Big\{C\subseteq B\mid f\in
Aut^{*}(\bar{N})\ \wedge\ f\restriction C=id_C\Rightarrow f\in
G\Big\}$$
\end{defi}

\begin{lemma}\label{TFAE}
Let $\delta$ be a limit ordinal and let $\bar{N}\in K_{\delta}^s.$ Let
$B\subseteq N_0$ 
such that $\big(N_0,c\big)_{c\,\in\, B}$ is saturated.
Let $C=acl\,C,\ C\subseteq B,$ and $g$ an elementary map with
$dom\,g=B,$ $g\restriction C=id_C,$
$\big(N_0,c\big)_{c\,\in\,B\,\cup\,g[B]}$ is saturated, and
$$B\,\indep{C}\,g(B)$$
Then the
following are equivalent. 

\begin{enumerate}
\item $C\in E$
\item All extensions of $g$ in $Aut^{*}(\bar{N})$ are in $G$
\item Some extension of $g$ in $Aut^{*}(\bar{N})$ is in $G$
\end{enumerate}

\end{lemma}
\proof $1.\Rightarrow 2.$ is trivial. 

\noindent $2.\Rightarrow 3.$ We just need to prove $g$ has some
extension in $Aut^{*}(\bar{N}).$ But this follows easily by the
saturation for every $j<\delta$ of $\big(N_{j+1},c\big)_{c\in N_j}.$ 

\noindent $3.\Rightarrow 1.$ Let $f\in Aut^{*}(\bar{N})$ such
that $f\restriction C=id_C.$ Let $n\in\omega$ and $g^*\in
Aut^*(\bar{N})$ such that $g^*\supseteq
g,\ f,g^*\in Aut\big(\bar{N}\restriction[n,\delta)\big),$ and
$g^*\in G.$ Let $B'\subseteq N_{n+1}$ such that
$B'\,\sindep{C}\,N_{n}$ and $tp(B'/C)=tp(B/C).$ Let $g_1\in
Aut\big(\bar{N}\restriction[n+2,\delta)\big)$ such that $g_1$
maps $g(B)$ onto $B'$ and $g_1\restriction B=id_{B}.$ Since
$g_1\restriction B=id_B,\ g_1\in G.$ Let $g_2=g_1g^*(g_1)^{-1}.$ Again
$g_2\in G,\ g_2\restriction C=id_C,$ and $g_2[B]=B'.$ As 
$$B'\,\indep{C}\,N_{n}$$
$f\in Aut(\bar{N}\restriction[n,\delta))$ and $f\restriction C=id_C,$
clearly
$$f(B')\,\indep{C}\,N_n$$
Therefore there exists $g_3\in
Aut\big(\bar{N}\restriction[n+2,\delta)\big)$ such that
$g_3\restriction B'=f\restriction B'$ and $g_3\restriction
N_{n}=id_{N_{n}},$ hence $g_3\in G.$  $(g_3)^{-1}f\restriction B'=id_{B'}$ so
$(g_2)^{-1}(g_3)^{-1}fg_2=id_{B}$ hence  $(g_2)^{-1}(g_3)^{-1}fg_2\in
G.$ But this implies $f\in G.$

\begin{theorem}\label{Shelah} 
Let $\big|T\big|<\lambda.$ Let $\bar{M}\in K_{\delta}^s.$ Let
$G\sub Aut^*(M).$ If
$$f\in Aut_{M_0}^*(\bar{M})\ \Rightarrow\ f\in G$$
but for no $C\sub M_0$ with $\big|C\big|<\lambda$  does
$$f\in Aut_{C}^*(\bar{M})\ \Rightarrow\ f\in G $$
then 
$$\big[Aut(M^*):G\big]>\lambda$$

\end{theorem}
\proof Suppose not. Let $\langle h_i\mid i<\lambda\rangle$ be a list
of the representatives of the left $G$ cosets of
$Aut(\bar{M}\restriction [1,\delta))$ possibly with repetition. Let
$\lambda=\bigcup\limits_{\zeta<cf\,\lambda}\lambda_{\zeta}$ with $\langle
\lambda_{\zeta}\mid \zeta<cf\,\lambda\rangle$ increasing continuous
and $\big|T\big|\leq\big|\lambda_0\big|\leq \big|\lambda_{\zeta}\big|
< \lambda.$ Let $M_0=\bigcup\limits_{\zeta< cf\,\lambda}M_{\zeta}^0$
and $M_1=\bigcup\limits_{\zeta< cf\,\lambda}M_{\zeta}^1$ with each
being a continuous chain such that $|M^i_{\zeta}|\leq
|\lambda_{\zeta}|.$

Now we define by induction on $\zeta<cf\,\lambda,$ $N_{0,\zeta},\
N_{1,\zeta},\ f_{\zeta},\
B_{\zeta},$ and $h_{j,\zeta}$ for $j<\lambda_{\zeta}$ such that

\begin{enumerate}
\item $f_{\zeta}$ is an automorphism of $N_{1,\zeta}$ 
\item $\langle f_{\zeta}\mid \zeta<cf\,\lambda\rangle$ is increasing
continuous
\item If $j<\lambda_{\zeta}$ and there is an $h\in
Aut(\bar{M}\restriction [1,\delta))$ such that
\begin{enumerate}
\item $h$ extends $f_{\zeta}$
\item $hG=h_jG$
\end{enumerate}
then $h_{j,\zeta}$ satisfies a. and b.
\item $B_{\zeta}$ is a subset of $N_{1,\zeta}$ of cardinality
$\leq\big|\lambda_{\zeta}\big|$
\item $M^1_{\zeta}\subseteq B_{\zeta}$
\item $N_{0,\zeta}\sub B_{\zeta+1}$ and $B_{\zeta+1}$ is closed
under $h_{j,\epsilon}$ and $h_{j,\epsilon}^{-1}$ for
$j<\lambda_{\epsilon}$ and $\epsilon\leq\zeta$
\item
$f_{\zeta+1}^{-1}(B_{\zeta+1})\,\sindep{N_{0,\zeta}}\,N_{0,\zeta+1}$
\item $N_{1,\zeta}\,\sindep{N_{0,\zeta}}\,M_0$
\item $M_1=\bigcup\limits_{\zeta< cf\,\lambda}N_{1,\zeta}\ \ \ \
M_0=\bigcup\limits_{\zeta< cf\,\lambda}N_{0,\zeta}$
\item $|N_{0,\zeta}|\leq |\lambda_{\zeta}|$
\item $(N_{1,\zeta+1},c)_{c\in N_{1,\zeta}}$ is saturted of
cardinality $\lambda$
\item $(M_1,c)_{c\in M_0\cup N_{1,\zeta}}$ is saturated of cardinality
$\lambda$
\end{enumerate}

\noindent For $\zeta=0$ let $B_0$ be empty, let $N_{0,0}$ be a
submodel of $M_0$ of cardinality $|\lambda_0|,$ let $N_{1,0}$ be a
saturated submodel of $M_1$ of cardinality $\lambda$ such that
$N_{1,0}\,\sindep{N_{0,0}}\,M_0$  and 
let $f_{\zeta}=id_{N_{1,0}}.$ At limit stages take unions. If
$\zeta=\epsilon+1,$ let $B_{\zeta}$ be as in 4,5,6. Let
$N_{0,\zeta}\subseteq M_0$ such that
$B_{\zeta}\,\sindep{N_{0,\zeta}}\,M_0,\ N_{0,\epsilon}\subseteq
N_{0,\zeta},\ M_{\zeta}^0\subseteq N_{0,\zeta},\ |N_{0,\zeta}|\leq
\lambda_{\zeta}.$ Let $N_{1,\zeta}\subseteq M_1$ such that
$B_{\zeta}\subseteq N_{1,\zeta},\
N_{1,\zeta}\,\sindep{N_{0,\zeta}}\,M_0,$ 
$(N_{1,\zeta},c)_{c\in N_{1,\epsilon}}$ is saturated of cardinality 
$\lambda,$ and
$(M_1,c)_{c\in M_0\cup N_{1,\zeta}}$ is saturated of cardinality 
$\lambda.$ Let $f_{\zeta}$
be an extension of $f_{\epsilon}\restriction N_{1,\epsilon}$ to an
automorphism of $N_{1,\zeta}$ so that
$$f_{\zeta}^{-1}(B_{\zeta})\,\indep{N_{1,\epsilon}}\,N_{0,\zeta}$$
Since
$$N_{0,\zeta}\,\indep{N_{0,\epsilon}}\,N_{1,\epsilon}$$
we have
$$f_{\zeta}^{-1}(B_{\zeta})\,\indep{N_{0,\epsilon}}\,N_{0,\zeta}$$
Let $f$ be an extension of
$\bigcup\limits_{\zeta<cf\,\lambda}f_{\zeta}$ to an element of
$Aut(\bar{M}\restriction [1,\delta)).$ We have defined
$f$ so that 

\begin{enumerate}
\item (By nonforking calculus) $\forall\zeta<cf\,\lambda,\ \forall j<\lambda_{\zeta},$
$$f^{-1}h_{j,\zeta}(M_0)\,\indep{N_{0,\zeta}}\,M_0$$
\item $f^{-1}h_{j,\zeta}\restriction N_{0,\zeta}=id$
\end{enumerate}

\noindent By lemma \ref{TFAE} none of the $f^{-1}h_{j,\zeta}$ are in
$G,$ a contradiction as for some $j<\lambda, \ fG=h_jG$ so for some
$\zeta, \ j<\lambda_{\zeta},$ $h_jG=h_{j,\zeta}G=fG.$

\begin{lemma}\label{niceNbar2}
Let $\big|T\big|<\lambda.$ Let $cf\,\delta\geq\kappa_r(T)+\aleph_1.$  
Suppose $\big[Aut(M^*):G\big]\leq\lambda$   and assume that for no $A\subseteq M^*$ with
$\big|A\big|<\lambda$ is $Aut_A(M^*)\subseteq G.$ Then for some
$\bar{N}\in K^s_{\delta},$
$$\bigwedge\limits_{\alpha<\delta}\,Aut_{N_{\alpha}}^{*}(\bar{N})\not\subseteq
G$$
\end{lemma}
\proof Suppose not. Let $\bar{M}\in K_{\delta}^s.$
Then there exists an $\alpha<\delta$ such that
$Aut^*_{M_{\alpha}}(\bar{M})\sub G.$ Without loss of generality $\alpha=0.$ 
By lemma \ref{Shelah} there exists
$E\subseteq  M_0$ such that $\big|E\big|<\lambda$ and
$Aut_E(\bar{M})\subseteq G.$ Let $f\in Aut_E(M^*)\backslash G.$ By
lemma \ref{build8} we can find $\bar{N}^1,\bar{N}^2\in
K^s_{\delta}$ and automorphisms $f_1\in
Aut_E(\bar{N}^1)$ and $f_2\in Aut_E(\bar{N}^2)$  such that

\begin{enumerate}
\item $E\subset N_0^1,\ E\subset N_0^2$
\item $f=f_2f_1$
\item $f_1\restriction E=f_2\restriction E=id_E$
\item $\forall\ \alpha,\beta<\delta,$  

\begin{enumerate}
\item $N_{\alpha}^1\,\sindep{N_{\alpha}^1\cap M_{\beta}}\,M_{\beta}$
\item $N_{\alpha}^2\,\sindep{N_{\alpha}^2\cap M_{\beta}}\,M_{\beta}$
\item $\big(N^1_{\alpha+1}\cap M_{\beta+1},\,c\big)_{c\in
(N^1_{\alpha+1}\cap M_{\beta})\cup (N_{\alpha}^1\cap M_{\beta+1})}$ is
saturated of cardinality $\lambda$
\item $\big(N^2_{\alpha+1}\cap M_{\beta+1},\,c\big)_{c\in
(N^2_{\alpha+1}\cap M_{\beta})\cup (N_{\alpha}^2\cap M_{\beta+1})}$ is
saturated of cardinality $\lambda$
\item $\big(N_{\alpha+1}^1\cap M_0\big)_{c\in N_{\alpha}^1\cap M_0}$
is saturated of cardinality $\lambda$
\item $\big(N_{\alpha+1}^2\cap M_0\big)_{c\in N_{\alpha}^2\cap M_0}$
is saturated of cardinality $\lambda$
\end{enumerate}

\end{enumerate}

Since $f\not\in G$ we can assume without loss of generality that
$f_1\not\in G.$ Also, by the hypothesis of suppose not we can assume
there is a $F\sub N^1_0$ such that $(N_0^1,c)_{c\in F}$ is saturated
and $Aut_F(\bar{N}^1)\subseteq G.$ By lemma \ref{Shelah} we can assume that
$\big|F\big|<\lambda$ and without loss of generality $E\subseteq F.$ Let
for $\alpha<\delta,$
$$F_{\alpha}=F\cap M_{\alpha}$$
By the lemma \ref{next}  we can find a sequence $\langle F_{\alpha}'\mid
\alpha<\delta\rangle$ such that for each $\alpha,$
$F_{\alpha}\subseteq F_{\alpha}'$ with $\big|F_{\alpha}'\big|<\lambda$
and for each $\beta<\alpha$ $F_{\alpha}'\cap M_{\beta}=F'_{\beta}$ and
if $F'=\bigcup\limits_{\alpha<\delta}F'_{\alpha}$ then 
$$M_{\alpha}\cap N_{0}^1\,\indep{F_{\alpha}'}\,F'$$
We define by induction on $\alpha<\delta$ a map
$g_{\alpha}$ an automorphism of $M_{\alpha}\cap N_{0}^1$ such that 

\begin{enumerate}
\item $\forall\,\beta,\alpha<\delta,\ \ \beta<\alpha\
\Rightarrow\ g_{\beta}\subseteq g_{\alpha}$
\item If $\alpha$ is a limit then
$g_{\alpha}=\bigcup\limits_{\beta<\alpha}g_{\beta}$
\item $g_{\alpha}(F'_{\alpha})\,\sindep{E}\,F'_{\alpha}$
\item $g_{\alpha}\restriction E=id_{E}$
\end{enumerate}

Let $\alpha=\beta+1$ and suppose $g_{\beta}$ has been defined. Let
$X\subseteq M_{\alpha}\cap N_{0}^1$ such that $X\frown
g_{\beta}(F'_{\beta})\equiv F'_{\alpha}\frown F'_{\beta}$ by
$h_{\beta}$ an extension of $g_{\beta}\restriction F_{\beta}'$ and
$$X\,\indep{g_{\beta}(F'_{\beta})}\,F'_{\alpha}\cup (M_{\beta}\cap N^1_{0})$$
Let $g_{\alpha}'=g_{\beta}\cup h_{\beta}.$ 
Since $X\,\sindep{g_{\beta}(F'_{\beta})}\,g_{\beta}(M_{\beta}\cap N_{0}^1)$ and
$F_{\alpha}'\,\sindep{F'_{\beta}}\,M_{\beta}\cap N_{0}^1,$ $g_{\alpha}'$
is an elementary map. Now let $g_{\alpha}$ be an extension of
$g_{\alpha}'$ to an automorphism of $M_{\alpha}\cap N_{0}^1.$ 
Let $g'=\bigcup\limits_{\alpha<\delta}g_{\alpha}.$ $g'$
is an automorphism of $N_0^1$ such that for every
$\alpha<\delta,$
$$g'[M_{\alpha}\cap N_{0}^1]=[M_{\alpha}\cap N_{0}^1]$$
By the saturation and independence of the $N_{\alpha}^1,\ M_{\beta}$
we can find  an extension $g$ of $g'$ such that $g\in Aut(\bar{N}_1)$
and $g\in Aut(\bar{M}).$ This gives a contradiction since
$g(F)\,\sindep{E}\,F$ and $g\in Aut(\bar{N}_1)$ implies $g\not\in G,$
but $g\in Aut(\bar{M})$ and $g\restriction E=id_E$ implies $g\in G.$

\begin{lemma}\label{next}
Let
$\bar{M}=\langle M_{\beta}\mid \beta\leq\delta\rangle\in
K^s_{\delta}.$  Let $F\subseteq M^*$
with $\big|F\big|<\lambda.$ Then there exits a set $F'$ such that
$\big|F'\big|<\lambda,\ F\subseteq F',$ and
$\forall\,\beta<\delta,$
$$*\ \ M_{\beta}\,\indep{F'\cap M_{\beta}}\,F'$$
\end{lemma}
\proof Let $w\subseteq F$ be finite. There are less than $\kappa_r(T)$
many $\alpha<\delta$ such that 
$$w\,\noindep{M_{\alpha}}\,M_{\alpha+1}$$
Let $a_w$ be the set of such $\alpha.$ 
For each $\alpha\in a_w$ let $w_{\alpha}\sub M_{\alpha}$ such that
$\big|w_{\alpha}\big|<\kappa_r(T),$ and
$$w\,\indep{w_{\alpha}}\,M_{\alpha}$$
Let $w^1=\bigcup\limits_{\alpha\in a_w} w_{\alpha}.$  Let
$F^1=\bigcup\limits_{w\mathop\subset\limits_{finite}F}w^1$ and repeat
this procedure $\omega$ times with $F^n$ relating to $F^{n+1}$ as $F$
is related to $F^1.$ Let $F'=\bigcup\limits_{n\in\omega}F^n.$ $F'$ satisfies $*.$

\begin{lemma}\label{3steps}
Let $Tr$ be a tree of infinite height. Let $\alpha< height(Tr)$ and let $\eta\in Tr\restriction
level(\alpha+1).$ Let $\langle M_{\beta}\mid \beta\leq \alpha\rangle$
be an increasing chain of models such that  for all $\beta<\alpha,$
$\big(M_{\beta+1},\ c\big)_{c\,\in\,M_{\beta}}$ is saturated. Let
$M_{\alpha}\subseteq N_0\subseteq N_1\subseteq N_2\subseteq N_3$ with
$\big(N_{i+1},\ c\big)_{c\,\in\,N_{i}}$ saturated for $i\leq 2.$
Suppose $\langle h_{\beta}\mid \beta\leq\alpha\rangle$ are such that

\begin{enumerate}
\item $h_{\beta}=id_{M_{\beta}}$
\item $h_{\beta}[N_i]=N_i$ for $i\leq 3$
\item $h_{\beta}[M_{\gamma}]=M_{\gamma}$ for $\gamma\leq\alpha$
\end{enumerate}

\noindent For each $\nu\in Tr\restriction level(\alpha+1)$ let
$m_{\nu},l_{\nu}$ be automorphisms of $N_0.$ Suppose $g_{\eta}\in
Aut(N_0)$ such that for all $\nu\in Tr\restriction level(\alpha+1),$
$$g_{\eta}m_{\eta}(m_{\nu})^{-1}(g_{\eta})^{-1}=l_{\eta}(l_{\nu})^{-1}h_{\gamma[\eta,\nu]}^{\eta(\gamma[\eta,\nu])\,<\,\nu(\gamma[\eta,\nu])}$$
Let $m_{\nu}^+,l_{\nu}^+$ be extensions of $m_{\nu}$ and $l_{\nu}$ to
automorphisms of $N_1$ for
all $\nu\in Tr\restriction level(\alpha+1).$ Then there exists a
$g'_{\eta}\in Aut(N_3)$ extending $g_{\eta}$ and for
all $\nu\in Tr\restriction level(\alpha+1)$ automorphisms of $N_3,$
$m_{\nu}'$ and $l_{\nu}'$ extending $m_{\nu}^+$ and $l_{\nu}^+$
respectively such that
$$g'_{\eta}m'_{\eta}(m'_{\nu})^{-1}(g'_{\eta})^{-1}=l'_{\eta}(l'_{\nu})^{-1}h_{\gamma[\eta,\nu]}^{\eta(\gamma[\eta,\nu])\,<\,\nu(\gamma[\eta,\nu])}$$
\end{lemma}
\proof Similar to the proof of lemma 1.8.

\begin{theorem}\label{last2}
Let $\big|T\big|<\lambda.$ Let $M^*$ be a saturated model of
cardinality $\lambda,$ and let 
$G\subseteq Aut(M^*).$  Suppose that for no $A\subseteq M$ with
$\big|A\big|<\lambda$ is $Aut_A(M^*)\subseteq G.$ Suppose $Tr$ is a tree
of height $\kappa,$ where $\kappa$ is a regular cardinal
$\geq\kappa_r(T)+\aleph_1$ such that each level of $Tr$ is of size at
most $\lambda,$ but $Tr$ having more than $\lambda$ branches. Then 
$$[Aut(M^*):G]>\lambda$$
\end{theorem}
\proof Suppose not. Then by lemma \ref{niceNbar2} there is a
$\bar{N}\in K_{\lambda\times\kappa}^s,$ such that 
$$\bigwedge\limits_{\alpha<\lambda\times
\kappa}Aut^*_{N_{\alpha}}(\bar{N})\not\subseteq G$$
By thinning $\bar{N}$ if necessary we can assume for each $\alpha<
\kappa$ there exists an automorphism $h_{\alpha}\in
Aut_{N_{\lambda\times\alpha}}(\bar{N})$ such that $h_{\alpha}\not\in
G.$ By induction  on $\alpha< \kappa$ for every $\eta\in
Tr\restriction level\,\alpha$ we define automorphisms
$g_{\eta},m_{\eta},l_{\eta}$ of $N_{\lambda\times\alpha}$ such that
if $\rho\neq\nu$ then $l_{\rho}\neq l_{\nu}$ and
$$g_{\rho}m_{\rho}(m_{\nu})^{-1}(g_{\rho})^{-1}=l_{\rho}(l_{\nu})^{-1}
h_{\gamma[\rho,\nu]}^{\rho(\gamma[\rho,\nu])\,<\,\nu(\gamma[\rho,\nu])}$$
At limit steps we take unions. If $\alpha=\beta+1,$ for each
$i<\lambda$ we define for some $\eta_i\in Tr\restriction
level\,\alpha,\ g_{\eta_i}\in Aut(N_{\lambda\times\beta+3i})$ such that
for each $\eta\in Tr\restriction level\,\alpha,\ \eta=\eta_i$ cofinally
many times in $\lambda,$ and for every $\nu\in Tr\restriction
level\,\alpha,$ $m_{\nu}^i\neq l_{\nu}^i\in Aut(N_{\lambda\times\beta
+3i})$ such that
$$g_{\eta_i}m_{\eta_i}^i(m_{\nu}^i)^{-1}(g_{\eta_i})^{-1}=l_{\eta_i}^i(l_{\nu}^i)^{-1} 
h_{\gamma[\eta_i,\nu]}^{\eta_i(\gamma[\eta_i,\nu])\,<\,\nu(\gamma[\eta_i,\nu])}$$
The $g_{\eta_i},\ m_{\nu}^i,\ l_{\nu}^i$ are easily defined by
induction on $i<\lambda$ using lemma \ref{3steps}. 
Then if we let
$g_{\eta}=\bigcup\big\{g_{\eta_i}\mid\eta_i=\eta\big\},\
m_{\eta}=\bigcup\limits_{i<\lambda}m_{\eta}^i$ and
$l_{\eta}=\bigcup\limits_{i<\lambda}l_{\eta}^i$ we have finished. Let
$Br$ the set of branches of $Tr$ of height $\kappa.$ For
$\rho\in Br$ let $g_{\rho}=\bigcup\big\{g_{\eta}\mid \eta<\rho\big\},\
m_{\rho}=\bigcup\big\{m_{\eta}\mid \eta<\rho\big\},$ and $ 
l_{\rho}=\bigcup\big\{l_{\eta}\mid \eta<\rho\big\}.$ If $\rho\neq\nu,$
$g_{\rho}\neq g_{\nu}$ since without loss of generality
$\rho(\gamma[\rho,\nu])<\nu(\gamma[\rho,\nu])$ and
$$g_{\rho}m_{\rho}(m_{\nu})^{-1}(g_{\rho})^{-1}=l_{\rho}(l_{\nu})^{-1}
h_{\gamma[\rho,\nu]}^{\rho(\gamma[\rho,\nu])\,<\,\nu(\gamma[\rho,\nu])}$$
and
$$g_{\nu}m_{\nu}(m_{\rho})^{-1}(g_{\nu})^{-1}=l_{\nu}(l_{\rho})^{-1}$$
implies
$$g_{\rho}(g_{\nu})^{-1}l_{\rho}(l_{\nu})^{-1}g_{\nu}(g_{\rho})^{-1}=
l_{\rho}(l_{\nu})^{-1}
h_{\gamma[\rho,\nu]}^{\rho(\gamma[\rho,\nu])\,<\,\nu(\gamma[\rho,\nu])}$$
So if $g_{\rho}=g_{\nu}$ this would imply
$h_{\gamma[\rho,\nu]}^{\rho(\gamma[\rho,\nu])\,<\,\nu(\gamma[\rho,\nu])}=id_{M^*}$
a contradiction. If
$$[Aut(M^*):G]\leq\lambda$$
then for some $\rho,\nu\in Br$ we must have $l_{\rho}(l_{\nu})^{-1}\in
G$ and $g_{\rho}(g_{\nu})^{-1}\in G,$ but then we get a contradiction
as 
$g_{\rho}(g_{\nu})^{-1}l_{\rho}(l_{\nu})^{-1}g_{\nu}(g_{\rho})^{-1}\in
G$ and $l_{\rho}(l_{\nu})^{-1}\in G,$ but
$h_{\gamma[\rho,\nu]}^{\rho(\gamma[\rho,\nu])\,<\,\nu(\gamma[\rho,\nu])}\not\in
G.$  

\begin{coro}
Let $G\subseteq Aut(M^*).$  Suppose that for no $A\subseteq M$ with
$\big|A\big|<\lambda$ is $Aut_A(M^*)\subseteq G.$ Suppose
$\big|T\big|<\lambda$ and $M^*$ does not
have the small index property. Then 

\begin{enumerate}
\item There is no tree of height an uncountable regular cardinal
$\kappa$ with at most $\lambda$ nodes, but more than $\lambda$
branches. 
\item For some strong limit cardinal $\mu,$ $cf\,\mu=\aleph_0$ and
$\mu<\lambda<2^{\mu}.$ 
\item $T$ is superstable.
\end{enumerate}

\end{coro}
\proof 
\begin{enumerate}
\item By the previous theorem
\item By 1. and [Sh 430, 6.3]
\item If $T$ is stable in $\lambda,$ then
$\lambda=\lambda^{<\kappa_r(T)},$ so if $\kappa_r(T)>\aleph_0$ we can
let $\kappa$ from the previous theorem be the least $\kappa$ such that
$\lambda<\lambda^{\kappa}.$
\end{enumerate}

\pagebreak

\begin{center}
REFERENCES
\end{center}

\begin{enumerate}
\item {[L] D. Lascar {\em The group of automorphisms of a relational
saturated structure,} (to appear)}

\item {[L Sh] D. Lascar and Saharon Shelah {\em Uncountable
Saturated Structures have the small index property,} (to appear)}

\item {[Sh 430] S. Shelah, {\em More cardinal arithmetic,} (to
appear)}

\item {[Sh c] S. Shelah, {\em Classification theory and the number
of isomorphic models, revised,} North Holland Publ. Co., Studies in
Logic and the foundations of Math, vol 92, 1990.}

\item  {[ShT] S. Shelah and S. Thomas, {\em Subgroups of small index
in infinite  symmetric groups II,} J Symbolic Logic 54 (1989)
1, 95-99}

\end{enumerate}

\end{document}